\documentclass[a4paper,12pt]{opelsarticle}
\usepackage{amsmath,bm,amsfonts,graphicx,tikz,pgfplots, pdfsync}
\usepackage[colorlinks=true, citecolor=red]{hyperref}
\usepackage{fullpage, graphicx}
\usetikzlibrary{patterns}

\usepackage{xcolor}
\usepackage{listings}

\definecolor{mGreen}{rgb}{0,0.6,0}
\definecolor{mGray}{rgb}{0.5,0.5,0.5}
\definecolor{mPurple}{rgb}{0.58,0,0.82}
\definecolor{backgroundColour}{rgb}{0.95,0.95,0.92}

\lstdefinestyle{CStyle}{
    backgroundcolor=\color{backgroundColour},   
    commentstyle=\color{mGreen},
    keywordstyle=\color{magenta},
    numberstyle=\tiny\color{mGray},
    stringstyle=\color{mPurple},
    basicstyle=\footnotesize,
    breakatwhitespace=false,         
    breaklines=true,                 
    captionpos=b,                    
    keepspaces=true,                 
    numbers=left,                    
    numbersep=5pt,                  
    showspaces=false,                
    showstringspaces=false,
    showtabs=false,                  
    tabsize=2,
    language=C
}

\newtheorem{definition}{Definition} 
\newtheorem{theorem}{Theorem} 
 
\newtheorem{proposition}{Proposition}

\newtheorem{remark}{Remark}
\newtheorem{conjecture}{Conjecture}

\def\beginproof{\par{\vspace{0.5em}\it Proof} \ignorespaces}
\def\endproof{{\hfill $\square$\\ \vspace{0.5em}}}

\def \R{{\mathbb{R}}}
\def \SS{{\mathbb{S}}}

\def \vae {{\varepsilon }}
\def \ds {\displaystyle}
\def\p{\partial}
\def\n{\nabla}

\def\d{\mathrm{d}}
\def\bkappa{\bar\kappa}
\def\vx {{\bm x}}

\def\omegav {{\bm\omega}}

\def\blangle {\Big\langle}
\def\brangle {\Big\rangle}
\def\oI {{\overline I}}

\begin{document}

\title{Remarks on the Radiative Transfer Equations for Climatology}
\author{Claude Bardos\footnote{University of Paris Denis Diderot (claude.bardos@gmail.com)},
Fran\c cois Golse\footnote{Ecole Polytechnique, Palaiseau 91128, France 
  ({francois.golse@polytechnique.edu}).}
and Olivier Pironneau\footnote{Applied Mathematics, Jacques-Louis Lions Lab, Sorbonne Universit\'e, 75252 Paris cedex 5, France
  ({olivier.pironneau@sorbonne-universite.fr}, \url{https://www.ljll.math.upmc.fr/pironneau/}).}}
  

\parindent=0pt
\begin{frontmatter}
\begin{abstract}
Using theoretical and numerical arguments we discuss some of the commonly accepted approximations for the radiative transfer equations in climatology.  
 \end{abstract}
\end{frontmatter}

\section{Introduction}
Satellite, atmospheric and terrestrial measurements are numerous to support the global warming.  Various hypothesis are made to explain the disastrous effect of the GreenHouse Gases (GHG).  However  a rigorous proof derived from the fundamental equations of physics is not available and one of the problems is the complexity of the physical system of planet Earth in its astronomical setting around the Sun.

In many references such as \cite{kaper},\cite{MOD},\cite{CHA} and \cite{FOW} it is explained that the Sun radiates light  with a heat flux $Q=1370$Watt/m$^{2}$,
in the frequency range $(0.5,20)\times 10^{14}$Hz corresponding approximately to a black body at temperature of 5800K; 70\% of this light intensity reaches the ground because the atmosphere is almost transparent to this spectrum and about 30\% is reflected back by the clouds or the ocean, snow, etc  (albedo). 
The Earth behaves almost like a black body at temperature $T_e=288K$ and as such radiates rays of frequencies $\nu$ in the infrared spectrum $(0.03,0.6)\times 10^{14}$Hz. 

``The absorption coefficient variation with frequencies is such that the atmosphere is essentially transparent to solar radiation'' (\cite{FOW}, p65).

``Solar radiation that is not absorbed or reflected by the atmosphere (for example by clouds) reaches the surface of the Earth. The Earth absorbs most of the energy reaching its surface, a small fraction is reflected\footnote{ \texttt{ https://public.wmo.int/en/sun\%E2\%80\%99s-impact-earth}}.''

Carbone dioxide renders the atmosphere opaque to infrared radiations  around $20$THz . Hence increasing its proportion in air increases the absorption coefficient in that range.  It is believed to be one of the causes of global warming.  However it must be more complex because  some experimental measurements show that increasing opacity decreases the temperature at high altitude, above the clouds \cite{DUF}.

In the framework of the radiative transfer equations, one can studies the  conjectures cited above.  To support theoretical claims, the numerical solver developed in \cite{CBOP}, \cite{OP} and \cite{FGOP} is used.  We will address 4 questions:
\begin{itemize}
\item What is the effect of increasing an altitude dependent absorption coefficient?
\item Is it true that sunlight crosses the Earth atmosphere unaffected?
\item Is  30\% of albedo equivalent to 30\% reduction of the Sun radiative energy ?
\item What is the effect on the temperature of increasing the absorption coefficient in a specific frequency range?
\end{itemize}

\section{The Radiative Transfer Equations for a Stratified Atmosphere}

Finding the temperature $T$ in a fluid heated by electromagnetic radiations is a complex problem because  interactions of  photons with  atoms in the medium involve  rather intricate quantum phenomena. 
Assuming local thermodynamic equilibrium leads to a well-defined electronic temperature. In that case, one can write a kinetic equation for the radiative intensity $I_\nu(\vx,\omegav,t)$ at time $t$, at position $\vx$ and
in the direction $\omega$ for photons of frequency $\nu$, in terms of the temperature field $T(\vx,t)$. 
\begin{equation}
\label{RTEq}
\begin{aligned}
\frac1c\p_t I_\nu + \omegav\cdot\n I_\nu+\rho\bkappa_\nu a_\nu\left[I_\nu-{\tfrac{1}{4\pi}\int_{\SS^2}} p(\omegav,\omegav')I_\nu(\omegav')\d\omega'\right]
\\
=\rho\bkappa_\nu(1-a_\nu) [B_\nu(T)-I_\nu]&.
\end{aligned}
\end{equation}
In this equation $c$ is the speed of light,  $\rho$ is the density of the fluid, $\n$ designates the gradient with respect to the position $\vx$, while
\begin{equation}\label{Planck}
B_\nu(T)=\frac{2 \hbar \nu^3}{c^2[{\rm e}^\frac{\hbar\nu}{k T}-1]}
\end{equation}
is the Planck function at temperature $T$, with $\hbar$ the Planck constant and $k$ the Boltzmann constant. Recall the Stefan-Boltzmann identity,
\begin{equation}\label{StefBoltz}
\int_0^\infty B_\nu(T)\d\nu=\bar\sigma T^4\,,\qquad\bar\sigma=\frac{2\pi^4k^4}{15 c^2\hbar^3}\,,
\end{equation}
where $\pi\bar\sigma$ is the Stefan-Boltzmann constant.

The intricacy of the interaction of photons with the atoms of the medium is contained in the mass-absorption $\bkappa_\nu$, which is theoretically the result of vibration and rotation atomic resonance, but for practical purpose the fraction of radiative intensity at frequency $\nu$ that is absorbed by fluid per unit length. 
The coefficient 
$a_\nu\in(0,1)$ is the scattering albedo, and $\frac1{4\pi}p(\omegav,\omegav')\d\omegav$ is the probability that an incident ray of light with direction $\omegav'$ scatters in the infinitesimal element of solid angle $\d\omegav$ centered at $\omegav$.

The kinetic equation \eqref{RTEq} is coupled to the temperature - or energy conservation - equations of the fluid .  When thermal diffusion is small, the system decouples and energy balance becomes:
\begin{equation}
\label{Energ}
\int_0^\infty\rho\bkappa_\nu(1-a_\nu)\left(\int_{\SS^2}I_\nu(\omegav)\d\omegav-4\pi B_\nu(T)\right)\d\nu=0\,.
\end{equation}
When the wave source is far in the direction $z$ the problem becomes one-dimensional in $z$ and the radiative intensity scattered in the direction $\omega$ depends only on $\mu$,  the cosine of the angle between $\omega$ and $Oz$. Consequently the system becomes,
\begin{equation}
\label{eqI}
\begin{aligned}
&\mu\p_z I_\nu + \rho\bkappa_\nu I_\nu 
= \rho\bkappa_\nu(1-a_\nu) B_\nu(T)+\frac{\rho\bkappa_\nu a_\nu}2\int_{-1}^1 p(\mu,\mu')I_\nu(z,\mu')\d\mu',
\\ & 
\int_0^\infty  \rho\bkappa_\nu(1-a_\nu)\left(B_\nu(T)-\tfrac12\int_{-1}^1 I_\nu\d\mu\right) \d\nu=0,\quad z\in(0,H),~~|\mu|<1,~~\nu\in\R^+\,.
\end{aligned}
\end{equation}
For this system to be  mathematically well posed, the radiation intensity $I_\nu$ must be given on the domain boundary where radiation enters. For example,
\begin{equation}
\label{oneabmu}
~I_\nu(H,-\mu)= Q^-(\mu), ~I(0,\mu)=Q^+(\mu),\quad 0<\mu<1\,.
\end{equation}

\subsection{Semi-analytical Solution when the scattering is isotropic}
Isotropic scattering is modelled by taking $p\equiv 1$. By introducing an optical depth (in \cite{FOW} and others, the signs are changed), 
\begin{equation}\label{optl}
\tau=\int_0^z\rho(\eta)d\eta,
\text{ $\rho$ becomes $1$  and $z$ becomes $\tau\in(0,Z)$, $Z=\int_0^{H}\rho(\eta)d\eta$.}
\end{equation}
Let us denote the exponential integral by
\begin{equation}\label{expintegr}
E_p(X):=\int_0^1e^{-X/\mu}\mu^{p-2}{d}\mu\,.
\end{equation}
Define 
\begin{equation}\label{JS}
\begin{aligned}
J_\nu(\tau) = \frac12\int_{-1}^1 I_\nu d\mu, 
\quad 
S_\nu(\tau) = \tfrac12\int_0^1\left(e^{-\frac{\kappa_\nu\tau}{\mu}}Q^+_\nu(\mu)+e^{-\frac{\kappa_\nu(Z-\tau)}{\mu}}Q^-_\nu(\mu)\right){d}\mu.
\end{aligned}
\end{equation}
The problem is equivalent to a functional integral equation:
\begin{equation}\label{eqint}
\left.
\begin{aligned}
&J_\nu(\tau)=S_\nu(\tau)+\frac{\kappa_\nu}2\int_0^ZE_1(\kappa_\nu\vert \tau-t\vert )\left(
(1-a_\nu)B_\nu(T(t)) + a_\nu J_\nu(t)\right)d t\,,
\\ 
&\int_0^\infty(1-a_\nu)\kappa_\nu B_\nu(T(\tau)){d}\nu=\int_0^\infty(1-a_\nu)\kappa_\nu  J_\nu(\tau){d}\nu,
\end{aligned}
\right\}
\end{equation}
Once $J_\nu, T$ are known, $I_\nu$ can be recovered by
\begin{equation}\label{IntForm}
\begin{aligned}
I_\nu(\tau,\mu)=&e^{-\frac{\kappa_\nu\tau}{\mu}}I_\nu(0,\mu)\mathbf 1_{\mu>0}+e^{\frac{\kappa_\nu(Z-\tau)}{\mu}}I_\nu(Z,\mu)\mathbf 1_{\mu<0}
\\
&+{\bf 1}_{\mu>0}\int_0^\tau e^{-\frac{\kappa_\nu|\tau-t|}{|\mu|}}\tfrac{\kappa_\nu}{\vert\mu\vert}\left(
(1-a_\nu)\kappa_\nu B_\nu(T(t)) + a_\nu\kappa_\nu J_\nu(t)\right)d t
\\ &
+{\bf 1}_{\mu<0}\int_\tau^Z e^{-\frac{\kappa_\nu|\tau-t|}{|\mu|}}\tfrac{\kappa_\nu}{\vert\mu\vert}\left(
(1-a_\nu)\kappa_\nu B_\nu(T(t)) + a_\nu\kappa_\nu J_\nu(t)\right)d t.
\end{aligned}
\end{equation}
\begin{remark}
Notice that two  solutions $(T,J_\nu),(T',J'_\nu)$ with $J_\nu\equiv J'_\nu$ will have  $T\equiv T'$ but not necessarily $I_\nu\equiv I_\nu'$.
\end{remark}

\section{Increasing an Altitude Dependent Absorption Coeficient Decreases the Temperature in the atmosphere}
The following argument goes against a common belief which claims that increasing the absorption coefficient - which GHG do indeed - implies an increase of temperature in the atmosphere. Consider the case of an altitude dependent absorption coefficient which we write as $r(z)\kappa$ with $\kappa$ constant.

\begin{proposition}\label{prop:1}
If $\kappa$ is independent of $\nu$ but function of the altitude  and $T$ is monotone decreasing with altitude, then, increasing $\kappa$ anywhere will decrease the temperature.
\end{proposition}

\beginproof~~
The proof is a straightforward consequence of the concept of optical length \eqref{optl}. 
Here, let 
$\tau=\int_0^z{r}(\zeta)d\zeta.
$
 Observe that $\partial_z=\frac{\partial\tau}{\partial \tau}\partial_\tau={r}(z)\partial_z$; consequently ${r}$ disappears from \eqref{eqI},
 $H$ is replaced by $Z=\int_0^H{r}(\zeta)d\zeta$, and $z$ is replaced by $\tau$.
 
When  $\kappa$ and  $a$ are constant,  $\ds\bar I_{r}=\int_0^\infty I_\nu d\nu$ is solution of the "grey model":
\begin{equation}
\label{eqIc}
\begin{aligned}
&\mu\partial_\tau \bar I_1 + \kappa\bar I_1 =\frac12\int_{-1}^1
\left(\kappa(1-a) \bar I_1  + \kappa a\int_{-1}^1 p(\mu,\mu'){\bar I_1}(\mu,\mu') d\mu' \right)d\mu, 
\\ &   I_{r}(H,-\mu)=0,\quad\bar I_{r}(0,\mu)=Q^+(\mu),\quad 0<\mu<1.
\end{aligned}
\end{equation}
where the Stefan-Boltzmann identity has been used:
$\ds\sigma T^4=\tfrac12\int_{-1}^1 {\bar I_{r}}\d\mu.
$

Let $\bar I_1(\tau,\mu)$ be the solution of \eqref{eqIc}.
We want to compare $I_{{r}+\delta{r}}$ with $I_{{r}}$ when $\delta{r}$ has a small support near altitude $y$. By \eqref{optl},
\begin{equation}
\label{eqId}
\begin{aligned}
&\bar I_{{r}+\delta{r}}(z) = \bar I_1(\int_0^z({r}(\zeta)+\delta{r}(\zeta))d\zeta)
\approx \bar I_1(\tau(z)) +  \frac{d}{d\tau}\bar I_1(\tau(y))\delta{r}(y)
\\ &
\implies\qquad
\frac12\int_{-1}^1\bar I_{{r}+\delta{r}}(z)d\mu - \frac12\int_{-1}^1\bar I_{r}(z) d\mu\approx \frac{d}{d\tau}\left(\frac12\int_{-1}^1\ \bar I_1(\tau) d\mu\right)|_{\tau(y)}\delta{r}(y)
\end{aligned}
\end{equation}
Finally from  \eqref{StefBoltz},
\begin{equation}
\label{eqIdT}
T_{{r}+\delta{r}}(z) - T_{{r}}(z) \approx 4\sigma T^3(y)\frac{\delta{r}(y)}{{r}(y)} \frac{d T_{r}}{d z}(y).
\end{equation}
\endproof

\section{Sunlight Crosses the Earth Atmosphere Unaffected?}
It is true that the absorption coefficient $\kappa_\nu$ is  small in the frequency range of
visible light; but does it imply that the light source can be used as a boundary condition at altitude zero instead of altitude H, the top of the troposphere?

First eliminate the $z$ dependency of the coefficients by using the optical length $\tau$ defined in \eqref{optl}.  
The following result gives an interesting answer to this question.

\begin{proposition}\label{prop:two}
Assume $a_\nu=0$ and assume that for some $\nu^*$, 
\[
\epsilon:=\max_{\nu>\nu^*}\kappa_\nu << 1 \text{ and }\quad \kappa_\nu Q^0_\nu|_{\nu>\nu^*}=O(1).
\]
Let $I_\nu,T$ and $I'_\nu,T'$ be solutions of  \eqref{eqI},\eqref{oneabmu} respectively with 
\begin{equation}
\begin{aligned}
& Q^+_\nu(\mu)=0,\quad Q^-_\nu(\mu)=\mu Q^0_\nu,\quad 0<\mu<1.
\\
& {Q'}^+_\nu=\mu {\rm e}^{-\frac{\kappa_\nu Z}\mu} Q^0_\nu, \quad {Q'}^-_\nu(\mu)=0,\quad 0<\mu<1.
\end{aligned}
\end{equation} 
\noindent  Then $ T'(\tau) = T(\tau) + O(\epsilon)$  and
\[
I'_\nu-I_\nu =  \mu {\rm e}^{-\frac{\kappa_\nu Z}\mu} Q^0_\nu\left({\bf 1}_{\mu>0} {\rm e}^{-\frac{\kappa_\nu}{\mu}\tau} + {\bf 1}_{\mu<0}{\rm e}^{\frac{\kappa_\nu}{\mu}\tau }\right) + O(\epsilon).
\]
\end{proposition}

\beginproof~~
Note that
\begin{equation}
\left.
\begin{aligned}\label{eqI2}
&\mu\partial_\tau  I_\nu + \kappa_\nu( I_\nu -B_\nu(T))=0, \quad\nu<\nu^*,
\\
 & I_\nu(0)|_{\mu>0}=0,\quad  I_\nu(Z)_{\mu<0}=Q^-_\nu(-\mu),\quad
 \quad\nu<\nu^*,
\\
&I_\nu\approx I_\nu^*(\mu)\text{ independent of $\tau$},~\quad \nu>\nu^*.
\\
& I^*_\nu|_{\mu>0}=0,\quad  I^*_\nu|_{\mu<0}=Q^-_\nu(-\mu),\quad
 \quad\nu>\nu^*.
\\
&\int_0^\infty \kappa_\nu B_\nu(T)d\nu
\approx\int_0^{\nu^*}\frac{\kappa_\nu}2\int_{-1}^1 I_\nu \,d\mu d\nu + \int_{\nu^*}^\infty\frac{\kappa_\nu}2\int_{-1}^0 Q^-_\nu(-\mu) \,d\mu d\nu
  \end{aligned}
  \right\}
\end{equation}

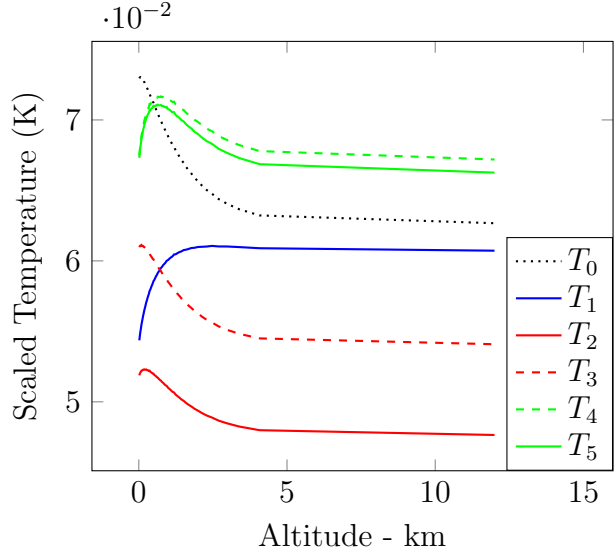
\begin{figure}[htbp]
\begin{center}
\begin{tikzpicture}[scale=1]
\begin{axis}[legend style={at={(1,0)},anchor=south east}, compat=1.3,
  xlabel= {Altitude - km},
  ylabel= {Scaled Temperature (K)},
  xmax=16]
%
\addplot[thick,dotted,color=black,mark=none,mark size=1pt] table [x index=0, y index=1]{fig/grey.txt};
\addlegendentry{ $T_0$}
\addplot[thick,solid,color=blue,mark=none,mark size=1pt] table [x index=0, y index=1]{fig/truethrough.txt};
\addlegendentry{ $T_1$}
\addplot[thick,solid,color=red,mark=none,mark size=1pt] table [x index=0, y index=1]{fig/correctedthrough.txt};
\addlegendentry{ $T_2$}
\addplot[thick,dashed,color=red,mark=none,mark size=1pt] table [x index=0, y index=1]{fig/through.txt};
\addlegendentry{ $T_3$}
\addplot[thick,dashed,color=green,mark=none,mark size=1pt] table [x index=0, y index=1]{fig/truethrough2.txt};
\addlegendentry{ $T_4$}
\addplot[thick,solid,color=green,mark=none,mark size=1pt] table [x index=0, y index=1]{fig/correctedthrough2.txt};
\addlegendentry{ $T_5$}
%
\end{axis}
\end{tikzpicture}
\end{center}
\caption{\label{figone}  Scaled temperatures  $z\to T(z)$ versus altitude, computed with $0.5(0.1+{\bf 1_{\nu<3}})$. The 3 curves $T_1,T_2,T_3$ corresponds to  {\bf 1:} $Q^-_\nu=Q^0_\nu$, $Q^+_\nu=0$, {\bf 2:} $Q^-_\nu=0$, $I_\nu(0,\mu)|_{\mu>0}= \mu{\rm e}^{-\frac{\kappa_\nu Z}\mu}Q^0_\nu$ and to {\bf 3:} $Q^-_\nu=0$, $Q^+_\nu=Q^0$.
These can be compared to the grey case $T_0$ with $\kappa=0.5$. Note that $T_1$ is very different from $T_2$. 
However if the Sun temperature is multiplied by 2 and $\kappa_\nu=0.5(0.1+{\bf 1_{\nu<0.2}})$, then the conditions in Proposition \ref{prop:two} are met and confirmed by the numerical simulations: $T_4$ is computed like $T_1$ and $T_5$ is computed like $T_2$. }
\end{figure}

Now  $I_\nu'$ is defined by \eqref{eqI} but with 
\[
I'_\nu(0)|_{\mu>0} = {Q'}_\nu^+(\mu),\qquad I'_\nu(Z)|_{\mu<0}=0.
\]
By the same argument, when $\nu>\nu^*$, $I'_\nu\approx {I'}^*_\nu(\mu)$, independent of $\tau$, and
\begin{equation}
\left.
\begin{aligned}\label{eqI3}
&\mu\partial_\tau  I'_\nu + \kappa_\nu( I'_\nu -B_\nu(T'))=0, \quad\nu<\nu^*,
\\
 & I'_\nu(0)|_{\mu>0} = {Q'}_\nu^+(\mu),\qquad I'_\nu(Z)|_{\mu<0}=	0,\quad
 \quad\nu<\nu^*,
\\
&\int_0^\infty \kappa_\nu B_\nu(T')d\nu
\approx\int_0^{\nu^*}\frac{\kappa_\nu}2\int_{-1}^1 I'_\nu d\mu d\nu 
+ \int_{\nu^*}^\infty\frac{\kappa_\nu}2\int_0^1 {Q'}^+_\nu (\mu)d\mu d\nu
\end{aligned}
\right\}
\end{equation}
Let $I^{\prime\prime}_\nu=I_\nu-I'_\nu$. It holds
\begin{equation}
\left.
\begin{aligned}\label{eqI4}
&\mu\partial_\tau  I^{\prime\prime}_\nu + \kappa_\nu( I^{\prime\prime}_\nu -(B_\nu(T)-B_\nu(T'))=0, \quad\quad\nu<\nu^*,
\\
&\int_0^\infty \kappa_\nu (B_\nu(T)-B_\nu(T'))d\nu=\int_0^{\nu^*}\frac{\kappa_\nu}2\int_{-1}^1 I^{\prime\prime}_\nu \,d\mu d\nu 
\\
&+ \int_{\nu^*}^\infty\left(\frac{\kappa_\nu}2\int_0^1 (Q^-_\nu(\mu)-{Q'}^+_\nu (\mu))\,d\mu\right)d\nu
\\
 & I^{\prime\prime}_\nu(0)|_{\mu>0} =-{Q'}_\nu^+(\mu),\qquad I^{\prime\prime}_\nu(Z)|_{\mu<0}=Q_\nu^-(-\mu),\quad
 \quad\nu<\nu^*.
  \end{aligned}
  \right\}
\end{equation}
We notice that $T=T'$ and $Q_\nu^-(\mu)={\rm e}^{-\frac{\kappa_\nu}\mu Z} {Q'}^+_\nu(\mu)$ implies
\begin{equation}
\left.
\begin{aligned}\label{eqI6}
&\mu\partial_\tau  I^{\prime\prime}_\nu + \kappa_\nu I^{\prime\prime}_\nu =0, \quad\nu<\nu^*,
\\
&
\int_0^{\nu^*}\left(\frac{\kappa_\nu}2\int_{-1}^1 I^{\prime\prime}_\nu \,d\mu\right)d\nu =
\int_{\nu^*}^\infty\left(\frac{\kappa_\nu}2\int_0^1 (1-{\rm e}^{-\frac{\kappa_\nu}\mu Z}){Q}^+_\nu (\mu)\,d\mu\right)d\nu
\\
 & I^{\prime\prime}_\nu(0)|_{\mu>0} = -{Q'}_\nu^+(\mu),\qquad I^{\prime\prime}_\nu(Z)|_{\mu<0}= {\rm e}^{\frac{\kappa_\nu}\mu Z}{Q'}_\nu^+(-\mu)
  \end{aligned}
  \right\}
\end{equation}
Let ${Q'}_\nu^+(\mu)=\mu Q^0_\nu$. Then 
\[
I^{\prime\prime}_\nu(0)|_{\mu>0} =-\mu Q^0_\nu,\qquad I^{\prime\prime}_\nu(Z)|_{\mu<0}=-\mu {\rm e}^{\frac{\kappa_\nu}\mu Z}Q^0_\nu .
\]
Consequently, from \eqref{eqI6} and \eqref{IntForm}, the solution is 
\begin{equation}
\begin{aligned}\label{eqI7}
&
I^{\prime\prime}_\nu(\tau,\mu) 
= - \mu Q^0_\nu\left({\bf 1}_{\mu>0} {\rm e}^{-\frac{\kappa_\nu}{\mu}\tau} 
+ {\bf 1}_{\mu<0}{\rm e}^{\frac{\kappa_\nu}{\mu}\tau }\right).
\\
\implies
&\int_{-1}^1 I^{\prime\prime}_\nu(\tau)d\mu =-Q^0_\nu E_3(\kappa_\nu \tau)+Q^0_\nu E_3(\kappa_\nu \tau)=0,
  \end{aligned}
\end{equation}
We have satisfied \eqref{eqI6} with precision $O(\kappa_\nu^2)$ because
\[
\kappa_\nu(1-{\rm e}^{-\frac{\kappa_\nu}\mu Z})\mu Q^0_\nu = \kappa_\nu^2 Z Q^0_\nu + O(\kappa_\nu^3)\approx 0, \quad \nu>\nu^*
\]
\endproof
\begin{remark}
It is easy to see that applying ${\rm e}^{-\frac{\kappa_\nu Z}\tau}\mu Q^0_\nu$ at $\tau=0$, leads
to change $E_3(\kappa_\nu\tau)$ in \eqref{IntForm} into  $E_3(\kappa_\nu(Z+\tau))$.
\end{remark}
The change in temperature is shown on Figure \ref{figone}. All variables are rescaled as in \cite{FGOPSIAM}$\S 7.1$. In particular temperatures are in Kelvin and divided by 4780. The air density is $1.22{\rm e}^{-z}$. Here $Q^0=5.66\cdot 10^{-5}$ and $T_{sun}=1.2$.  

\emph{ Conclusion: For the specific case of sunlight on Earth It is a gross approximation to assume that sunlight crosses the atmosphere unaffected.}  

However the theory works when $T_{sun}=2.4$ and $\nu^*=0.2$, probably because $\max_\nu B_\nu(T_{sun}\approx 18$ against 2 with $T_{sun}=1.2$.

\section{Is the 30\% Earth Albedo Equivalent to a 30\% reduction of the Sun Radiative Energy?}

Note that \eqref{IntForm} can be used to find the solution of \eqref{eqI}-a with generalized albedo (inspired by the model proposed in \cite{DUF}):
\[
I_\nu(0,\mu) = \sum_k \alpha_k I_\nu(\tau_k,-\mu) + Q^0_\nu\mu, \quad I_\nu(Z,-\mu)=0,\quad 0<\mu<1.
\]
with $k=1,2,..K$.  Indeed by \eqref{IntForm}
\begin{equation}\label{IntForm3}
\begin{aligned}
I_\nu(\tau,\mu)=&e^{-\frac{\kappa_\nu\tau}{\mu}}\left(\sum_k\alpha_k I_\nu(\tau_k,-\mu) + Q^0_\nu\mu\right)\mathbf 1_{\mu>0}
\\
&+{\bf 1}_{\mu>0}\int_0^\tau e^{-\frac{\kappa_\nu|\tau-t|}{|\mu|}}\tfrac{\kappa_\nu}{\vert\mu\vert}\left(
(1-a_\nu)\kappa_\nu B_\nu(T(t)) + a_\nu\kappa_\nu J_\nu(t)\right)d t
\\ &
+{\bf 1}_{\mu<0}\int_\tau^Z e^{-\frac{\kappa_\nu|\tau-t|}{|\mu|}}\tfrac{\kappa_\nu}{\vert\mu\vert}\left(
(1-a_\nu)\kappa_\nu B_\nu(T(t)) + a_\nu\kappa_\nu J_\nu(t)\right)d t
\end{aligned}
\end{equation}
\begin{equation}
\begin{aligned}
&=\mathbf 1_{\mu>0} e^{-\frac{\kappa_\nu\tau}{\mu}}\sum_k\alpha_k
\int_{\tau_k}^Z e^{-\frac{\kappa_\nu(t-\tau_k)}{\mu}}\tfrac{\kappa_\nu}{\mu}\left(
(1-a_\nu)\kappa_\nu B_\nu(T(t)) + a_\nu\kappa_\nu J_\nu(t)\right)d t
\\
&+{\bf 1}_{\mu>0}\left(e^{-\frac{\kappa_\nu\tau}{\mu}}Q^0_\nu\mu +\int_0^\tau e^{-\frac{\kappa_\nu|\tau-t|}{|\mu|}}\tfrac{\kappa_\nu}{\vert\mu\vert}\left(
(1-a_\nu)\kappa_\nu B_\nu(T(t)) + a_\nu\kappa_\nu J_\nu(t)\right)d t\right)
\\ &
+{\bf 1}_{\mu<0}\int_\tau^Z e^{-\frac{\kappa_\nu|\tau-t|}{|\mu|}}\tfrac{\kappa_\nu}{\vert\mu\vert}\left(
(1-a_\nu)\kappa_\nu B_\nu(T(t)) + a_\nu\kappa_\nu J_\nu(t)\right)d t
\end{aligned}
\end{equation}
\begin{equation*}
\begin{aligned}
\implies
\\&
J_\nu(\tau) = \tfrac12 Q^0_\nu E_3(\kappa_\nu\tau) 
+\tfrac12\int_0^ZE_1(\kappa_\nu\vert \tau-t\vert )\left[
(1-a_\nu) B_\nu(T(t)) + a_\nu\kappa_\nu J_\nu(t)\right]d t\,,
\\& 
+ \tfrac12\sum_k\alpha_k \kappa_\nu\int_{\tau_k}^Z E_1(\kappa_\nu(t+\tau-\tau_k))
\left[
(1-a_\nu)\kappa_\nu B_\nu(T(t)) + a_\nu\kappa_\nu J_\nu(t)\right]d t.
\end{aligned}
\end{equation*}
\subsection{Earth Albedo}
Consider the case
\[
 I_\nu(Z,\mu)|_{\mu<0}=0,\quad I_\nu(0,\mu)|_{\mu>0}=Q^0_\nu + \alpha I_\nu(0,-\mu)\,.
 \]
Assume no scattering, $\tau_1=0$, $K=1$.
As in \cite{FGOP}, the following iterative scheme is considered:
\begin{equation}\label{inteq}
\begin{aligned}
J_\nu^{n+1}(\tau) &= 
 \tfrac12 Q^0_\nu E_3(\kappa_\nu\tau) 
+\frac{\kappa_\nu}2\int_0^Z \left(E_1(\kappa_\nu\vert \tau-t\vert )+\alpha  E_1(\kappa_\nu(t+\tau))\right)B_\nu(T^n(t)) d t\,,
\\ &
\ds\int_0^\infty \kappa_\nu B_\nu({T^{n+1}_\tau})\d\nu = \int_0^\infty \kappa_\nu J^{n+1}_\nu(\tau)\d\nu\,.
\end{aligned}
\end{equation}

\subsection{Numerical Example}

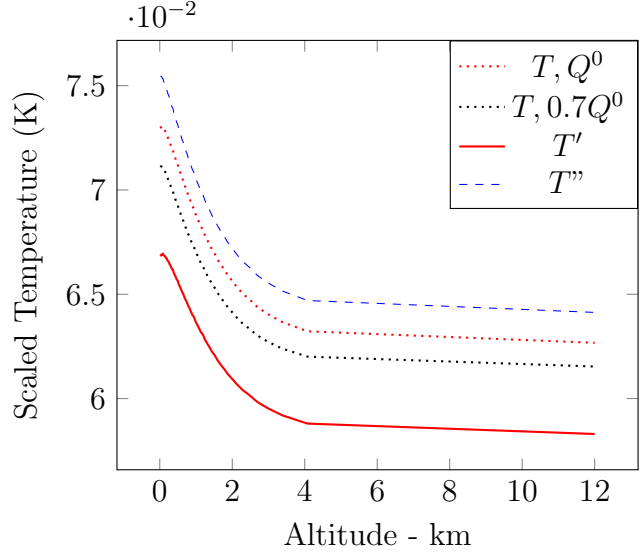
\begin{figure}
\begin{center}
\begin{tikzpicture}[scale=1]
\begin{axis}[legend style={at={(1,1)},anchor=north east}, compat=1.3,
  xlabel= {Altitude - km},
  ylabel= {Scaled Temperature (K)}]
%
\addplot[thick,dotted,color=red,mark=none,mark size=1pt] table [x index=0, y index=1]{fig/grey.txt};
\addlegendentry{ $T, Q^0$}
\addplot[thick,dotted,color=black,mark=none,mark size=1pt] table [x index=0, y index=1]{fig/Q07.txt};
\addlegendentry{ $T,0.7Q^0$}
\addplot[thick,solid,color=red,mark=none,mark size=1pt] table [x index=0, y index=1]{fig/mirror3.txt};
\addlegendentry{ $T'$}
\addplot[dashed,color=blue,mark=none,mark size=1pt] table [x index=0, y index=1]{fig/albedo03Q0.txt};
\addlegendentry{ $T"$}
%
\end{axis}
\end{tikzpicture}
\end{center}
\caption{\label{figtwo}  Scaled temperatures  $z\to T(z)$ versus altitude.  In all cases $Q^-_\nu=0$. The doted curves are $T$ when  $Q^+_\nu = \lambda Q^0_\nu$ with $\lambda=1$ (top dotted curve) and $\lambda=0.7$ (lower dotted curve); the doted line is the solution $T'$ when   $I_\nu(0,\mu)|_{\mu>0}= 0.3 I_\nu(0,-\mu)+\mu{\rm e}^{-\frac{\kappa_\nu Z}\mu}Q^0_\nu/0.7$ .  The dashed line is  when  $I_\nu(0,\mu)|_{\mu>0}= 0.3 I_\nu(0,-\mu)+\mu Q^0_\nu$.}
\end{figure}

Assume  $a_\nu=0$.
Figure \ref{figtwo} compares the numerical  solutions of \eqref{eqI} with $\kappa_\nu=0.5[{\bf 1}_{\nu<6} +0.1]$ and
\[
Q^-_\nu=0,\quad Q^+_\nu=\lambda Q^0_\nu,\quad \lambda=1\text{ or }0.7.
\]
with the solution of \eqref{eqI} with $\beta=0$ or 1 in
\[
{{Q^-_\nu}}^\prime=0,\quad {{I_\nu}}^\prime(0,\mu) = 0.3 {{I_\nu}}^\prime(0,-\mu) + \mu  Q^0_\nu[\tfrac 1{0.7}{\rm e}^{-\frac{\kappa_\nu Z}\mu}\beta+ 1-\beta],\quad 0<\mu<1.
\]

This numerical simulation is an attempt to replace the radiative source at $\tau=H$ by a radiative source at $\tau=0$ which gives the same result.  

Notice alongside that adding 0.3 albedo induces a heating of the atmosphere (i.e. comparing both computations with $\lambda=1,\beta=1$ above).

\section{Earth Albedo, General Statement, Accommodation Coefficient} 
As above one considers a stratified atmosphere \eqref{eqI} with $(z,\mu,\nu) $ in $(0,H)\times(-1,1) \times \R^+$; the radiation intensity $I_\nu(z,\mu)$  satisfies on at $z=H$  an incoming condition 
\begin{equation}
 I_\nu(H,-\mu)= Q_\nu^-(\mu), \quad 0<\mu <1,
\end{equation}
and on Earth, at $z=0$, some ``albedo" condition, denoted:
\begin{equation}
\begin{aligned}
I_{\nu,}(0,\mu)=\mathcal A_\nu(I_{\nu}(0,-\mu ),\quad 0<\mu <1,
\end{aligned}
\end{equation}
with $\mathcal A$ representing the albedo effect of the earth, hence being an operator from the space of outgoing intensities  into the space of incoming intensities.  

Below will be given, based on energy  estimates, some sufficient conditions on $\mathcal A$ for existence, stability -- hence uniqueness -- of the solution, beginning with the grey model and then considering a natural  operator which involves the effective temperature of Earth (see \cite{FOW} page  66).

\subsection{The energy type estimate for the grey problem}

For clarity we present the case without scattering ($a_\nu=0$) and with a constant $\kappa_\nu$ independent of the frequency $\nu$. By rescaling $z$ with $\kappa $,  the problem 
can be formulated in term of
$\ds 
\overline I(z,\mu) = \int_0^\infty I_\nu(z,\mu,) \d\nu
$
as, 
\begin{subequations}
\begin{equation}\label{albedo.Milne}
 \mu\partial_z \overline I(z,\mu) +\overline I(z,\mu)-\frac12\int_{-1}^1 \overline I(z,\mu')\d \mu'=0\,,
\end{equation}
\begin{equation}\label{albedo.Milne.bd}
 \overline I(H,-\mu) =Q(\mu):=\overline Q^-_\nu(\mu), \quad \overline I(0,\mu)=\mathcal A ( \oI(0,-\mu)), \quad 0<\mu<1\,.
\end{equation}
\end{subequations}
Note that we have assumed that ${\mathcal A}$ commutes with the integration in $\nu$.
First observe that the operator 
$\ds
\oI \mapsto \oI -\frac12\int_{-1}^1 \oI (\mu') \d \mu' 
$
is, in the space $L^2(-1,1)$ the orthogonal projection on functions of mean value $0$, because
\begin{equation}\label{pos}
\int_{-1}^1( \oI -\frac12\int_{-1}^1 \oI (\mu') \d \mu') \oI(\mu) d\mu = \int_{-1}^1(\oI -\frac12\int_{-1}^1 \oI (\mu') \d \mu') ^2  \d\mu \,.
\end{equation}
This implies that the left hand side of (\ref{pos}) is nonnegative and equal to $0$ iff 
$\ds \oI = \frac12\int_{-1}^1 \oI (\mu') \d \mu' 
$
is independent of $\mu\,.$
Multiplying the equation (\ref{albedo.Milne}) by $\oI(z,\mu) \,,$ integrating over $(0,H)\times(-1,1)\,,$ using Green's formula and inserting in the computations the conditions (\ref{albedo.Milne.bd}) one obtains:
\begin{equation} \label{pos2}
\begin{aligned}
&\int_0^H \d z\int_{-1}^1( \oI -\frac12\int_{-1}^1 \oI (\mu') \d \mu') ^2  \d\mu +\frac12 \int_0^1 \mu ( \oI(H,-\mu)^2 \d \mu 
\\
&+ \frac12 \int_0^1 \mu\big( ( \oI(0,\mu))^2-  (\mathcal A(\oI(0,\mu))^2 \big) \d \mu 
 =\frac12 \int_0^1 \mu (Q(\mu))^2 \d \mu 
\end{aligned}
\end{equation}
This estimate indicates that any property which implies the positivity of 
\[
\int_0^1 \mu\big( \oI(0,\mu)^2-  (\mathcal A(\oI (\mu, 0))^2 \big) \d \mu 
\] 
will lead to a well posed problem with a unique nonnegaative solution;  for instance:

\begin{theorem} \label{hilbert} 
Assume that the operator $\oI\mapsto \mathcal A(\oI)$ restricted to the positive cone $C^+(\oI(\mu)\ge 0)$ on the space $(L^2(-1,1),\mu,d\mu)$ then for any   positive $\mu\mapsto Q(\mu)$ there is a unique nonnegative solution of the
(\ref{albedo.Milne}),(\ref{albedo.Milne.bd}).
\end{theorem}
\begin{remark}\label{infty} Making use of the maximum principle one could prove similar results when the albedo operator is a contraction in the positive cone of  $L^\infty(0,1)\,.$ 
\end{remark}
\begin{remark}
The albedo condition which is a relation, for $0<\mu<1\,,$ between $\oI(0,\mu)\big|_{\mu>0}$ 
and $\oI(0,\mu)|_{\mu<0}$ can also be expressed in term of $\oI(0,\mu)|_{\mu>0}$ only.
 \end{remark}
Indeed,
define the operator ${\mathcal T}: \oI_+(\mu)\mapsto {\mathcal T}(\oI_+(\mu)) = I(0,-\mu)|_{\mu>0}$ (defined for instance in $L^\infty(0,1))$ by solving (\ref{albedo.Milne}) with
\begin{equation}\label{albedo.Milne.bdd}
 \oI(H,-\mu) =Q(\mu),  \quad \oI(0,\mu) =\oI_+(\mu),\quad 0<\mu<1  \,.
 \end{equation}
 Then one has:
 \begin{equation}
\oI(0,\mu)|_{\mu>0}=\mathcal A (\oI(0,-\mu))\Leftrightarrow \oI(0,\mu)|_{\mu>0}=\mathcal A ( \mathcal T(\oI(0,\mu)|_{\mu>0}))\,.
 \end{equation}
 With $Q(\mu)>0$ in $L^\infty(0,1)$, $\mathcal T$ is an affine operator which preserves the positivity and the monotonicity. 
 and this leads to the simple but  useful proposition, which , in turn, is the key to the proof of Theorem \ref{hilbert}.
 \begin{proposition}
 Consider two solutions $\{I_i\}_{i=1,2}$ of the (\ref{albedo.Milne})  with (\ref {albedo.Milne.bd}) with the same $Q(\mu)$ but  with two different albedo operators $\mathcal A_i\,.$ 
 Assume that both $\mathcal A_i$ are   linear contractions and  one of them is a strict linear contraction, which preserve positivity.  If 
 \begin{equation}\label{order}
 \forall   f\ge 0 \in L^\infty(0,1) \quad \forall \mu \in (0,1)\,, \quad   A_2(f)(\mu)\le A_1(f)(\mu)
 \end{equation}
 then one has the same ordering for the corresponding solutions :
 \begin{equation}
 \forall (z,\mu) \in (0,1)\times (0,H)\quad I_1(z,\mu)\le I_2(z,\mu)\,.
 \end{equation}
 \end{proposition}
\beginproof~~ 
 Using the linearity of the operators $\mathcal A_i$ one observes that $R =I_1-I_2$ solves\eqref{albedo.Milne} with $R(H,-\mu)=0$ and, for all $0<\mu<1$,
\begin{equation}
R(0,\mu)-\frac12(\mathcal A_1+ \mathcal A_2)(\mathcal T( R(0,\mu)) =   \frac12(\mathcal A_1- \mathcal A_2)(\mathcal T ((I_1(0,\mu)+I_2(0,\mu)))\,.
\end{equation}
Moreover, restricted to   solutions which satisfy $ R(H,-\mu) =0$, $\mu>0$ the operator $\mathcal T$ is a strict monotonicity preserving linear contraction; the same observation holds for 
$R\mapsto \frac12(\mathcal A_1+ \mathcal A_2)(\mathcal T (R)$.
As a consequence   the operator 
$R\mapsto (I-  \frac12(\mathcal A_1+ \mathcal A_2)(\mathcal T))$ is invertible with inverse given by the Neumann series (also preserving the positivity):
\begin{equation}\label{order2}
R= \sum_{k\ge 0} (\frac12(\mathcal A_1+ \mathcal A_2)(\mathcal T)))^k \frac12(\mathcal A_1- \mathcal A_2)(\mathcal T ((I_1+ I_2)|_{z=0,\mu>0}\,. 
\end{equation} 
$I_1$ and $I_2$ are positive intensities and for 
$\frac12(\mathcal A_1- \mathcal A_2)(\mathcal T ((I_1+ I_2)|_{z=0,\mu>0}$
and  the righthandside of (\ref{order2}), (\ref{order}) is also true.
\endproof

\subsection{The frequency dependent case}
 
Let us come back to \eqref{eqI} with \eqref{albedo.Milne.bd}. Let us denote $f_+=\max(f,0)$. When not ambiguous, let us write $I_\nu(\mu)$ instead of $I_\nu(0,\mu)$.
\begin{definition}
{\it An operator $\mathcal A$ defined in $L^1_\mu(0,1)$ is non-accretive if}
\begin{equation}\label{weakaccretive}
\forall I^1,I^2  \in  L^1_\mu(0,1), ~\int_0^1 \mu\left((I^2-I^1)_+-(\mathcal A(I^2)-\mathcal A(I^1))_+ \right)\d \mu \ge 0\,.\end{equation}
\end{definition}
\begin{remark}
Obviously (\ref{weakaccretive}) follows from the same pointwise property:
\begin{equation}
\forall (I^1,I^2) \quad ((\mathcal A(I^2)-\mathcal A(I^1))_+ \leq (I^2-I^1)_+,\quad \mu>0\,.
\end{equation}
\end{remark}
Let us uses the symbol 
$\blangle  f\brangle=\int_0^\infty \tfrac12 \int_{-1} ^1 f_\nu (\mu) \d \mu \d \nu
$.

 Then for two solutions $(I^2_\nu, T^2)$, $(I^1_\nu, T^1)$, by \eqref{eqI}, according to Theorem 4.1 in \cite{FGOP} one has:
 \begin{equation}
\begin{aligned}
\begin{aligned}
-\frac{d}{dz}\langle\mu(I^2_\nu-I^1_\nu)_+\rangle=
\\
=\left\langle\rho\bar\kappa_\nu(1-a_\nu)((I^2_\nu-I^1_\nu)-(B_\nu(T^2)-B_\nu(T^1)))\mathbf 1_{I^2_\nu>I^1_\nu}\right\rangle
\\
+\left\langle\rho\bar\kappa_\nu a_\nu\left((I^2_\nu-I^1_\nu)-\tfrac12\int_{-1}^1p(\mu,\mu')(I^2_\nu-I^1_\nu)(z,\mu')d\mu'\right)\mathbf 1_{I^2_\nu>I^1_\nu}\right\rangle
\ge 0
\end{aligned}
\end{aligned}
\end{equation}
and therefore 
$
\frac{d}{dz} \blangle  \mu ( I^2_\nu(z,\mu) - I^1_\nu(z,\mu) )_+\brangle  \leq 0\,.
$ 
On the other hand 
\[
I_\nu^2(H,-\mu) \le I_\nu^1(H,-\mu), \quad \mu>0\quad\implies~~
\blangle \mu(I_\nu^2(H,\mu)-I^1_\nu(H,\mu))_+\brangle = 0
\]
while the non accretivity of $\mathcal A $ gives 
\begin{equation}\label{firstmoment}
\begin{aligned}
\begin{aligned}
\langle\mu(I^2_\nu(0,\mu)-I^1_\nu(0,\mu))_+\rangle
\\
=\tfrac12\int_0^\infty d\nu\int_0^1\mu(\mathcal A(I^2_\nu(0,-\mu))-\mathcal A(I^1(0,-\mu)))_+-(I^2_\nu-I^1_\nu)_+(0,-\mu))d\mu\le 0
\end{aligned}
\end{aligned}
\end{equation}
Therefore, unless $\langle\mu(I^2_\nu(0,\mu)-I^1_\nu(0,\mu))_+\rangle=0$ there is a contradiction because it is a function which is decreasing, negative at $z=0$ and zero at $z=H$.
Then consider the difference of  equations (\ref{eqI}) with arguments 
$(I_\nu^i, T^i)\, i=1, 2$ multiplied by $ 
\frac{\mu}{\bar\kappa_\nu}
{\bf 1}_{I^2_\nu>I^1_\nu}
$.  Following the end of the proof of Theorem 4.1 in \cite{FGOP} (see also \cite{G87}), one finds by the same method  that existence and uniqueness derives from
$$
\frac{d}{dz}\left\langle\frac{\mu^2}{\bar\kappa_\nu}(I^2_\nu-I^1_\nu)_+\right\rangle=0\,,\qquad 0<z<H\,.
$$
\endproof
Eventually one has: 

\begin{theorem} \label{nudep}
Consider \eqref{eqI} with  incoming data $0\le I_\nu(H,-\mu)\le B_\nu (T_M)$, depending on a given (but possibly large) temperature $T_M$, and with outgoing data given by an accretive operator, of the form $I_{\nu}(0,\mu)|_{\mu>0}={\mathcal A} (I_{\nu}(0,-\mu))|_{\mu>0}$. The problem
has a unique well defined nonnegative solution.
\end{theorem}


\subsection{Examples of Albedo operators}
Below are given some examples of albedo operators which combine an accommodation parameter $\alpha$ and the reflection and the thermalisation effects of the Earth.

The simplest example  would be given (cf. Remark \ref{infty}) with  $0\le \alpha \le $ by the formula, $p>0$:
\begin{equation}
\mathcal A(I(0,\mu))  = \alpha I(0,-\mu) + (1-\alpha) \int_0^1 \mu^p I(0,-\mu') \d \mu'\,.
\end{equation}
For the  more realistic frequency dependent  cases, one should compare the reflection with the emission from the Earth as global black body under an effective temperature 
  $T_e\simeq  288$  and this leads to try the albedo operator
\begin{equation}
 \mathcal A(I_\nu(0,-\mu)) = \alpha I_\nu(0,-\mu) 
 +(1-\alpha) B_\nu(T_e) ,~~\mu>0. 
\end{equation}
This obviously satisfies the hypothesis of non accretivity.   

\section{ Calculus of Variations for the $\nu$-dependent Case}
The following argument sheds some light on the conditions needed to obtain a cooling (resp. heating)  from a local increase of the absorption coefficient.

\begin{conjecture}
Let $I^\vae_\nu(\tau,\mu), T^\vae(\tau)$ be the solution of  \eqref{eqI} when 
\begin{equation}\label{dkappa}
\kappa^\epsilon_\nu:=\kappa+\epsilon\delta\kappa_\nu\,,
\end{equation} 
Let $I^0_\nu(\tau,\mu), T^0(\tau)$ be the solution of  \eqref{eqI} with $\kappa_\nu=\kappa$ constant.
When $\delta\kappa_\nu\approx { \bm \delta_{\nu^*}}  $, the Dirac mass at $\nu^*$, and $\kappa$ is small, then the sign of $\ds \frac{ d T}{d\vae}(\tau)|_{\vae=0}$ is governed by the sign of 
\begin{equation}\label{sign}
\frac12\int_{-1}^1 I^0_{\nu^*}(\tau)d\mu - B_{\nu^*}(T^0(\tau)).
\end{equation}
\end{conjecture}
\subsection{Calculus of Variations Support of the Conjecture}

It will be convenient to define the Planck function in terms of the quantity given by the Stefan-Boltzmann law
$
\Phi:=\sigma T^4\,.
$
Henceforth we denote
$$
B_\nu(T)=b_\nu(\Phi)\,.
$$

We seek to study how the average radiation intensity or the temperature is altered if the absorption is modified on various intervals in the frequency variable. 
We expect that the simplest situation corresponds to absorption of the form \eqref{dkappa} with $\kappa>0$ independent of $\nu$, while $0<\epsilon\ll 1$. This is of course an extremely general formulation, but in practice one could think of
$$
\delta\kappa_\nu:=\kappa\mathbf 1_{\nu_1<\nu<\nu_2}\,.
$$
This corresponds to multiplying the absorption by $(1+\epsilon)$ in the frequency interval $(\nu_1,\nu_2)$ only, and leaving it invariant for all the other frequencies.

We assume that the scattering is isotropic with constant rate $\lambda:=\kappa a\ge 0$.

Henceforth, we consider the radiative transfer equation on $(0,Z)\times(-1,1)\times\R$,
\begin{equation}\label{RTeps}
\mu\partial_\tau I^\epsilon_\nu+\kappa^\epsilon_\nu(I^\epsilon_\nu-b_\nu(\Phi^\epsilon))+\lambda(I^\epsilon_\nu-\widetilde I^\epsilon_\nu)=0\,,\qquad
\blangle\kappa^\epsilon_\nu(\widetilde I^\epsilon_\nu-b_\nu(\Phi^\epsilon))\brangle=0\,,
\end{equation}
with the boundary conditions
\begin{equation}
I^\epsilon_\nu(0,\mu)=\mu Q^0_\nu,\quad I^\epsilon_\nu(Z,-\mu)=0\,,\qquad\qquad 0<\!\mu\!<1,
\end{equation}
and the notations
$$
\widetilde\phi:=\tfrac12\int_{-1}^1\phi(\mu)d\mu\,,\qquad\blangle\psi\brangle:=\int_0^\infty\psi(\nu)d\nu\,.
$$

First we study the case $\epsilon=0$: since $\kappa$ is constant, one can average in frequency, and set 
$$
{\bar I^0}(\tau,\mu):=\blangle I^0_\nu(\tau,\mu)\brangle
,
 \qquad\widetilde {\bar I^0}(\tau)=\blangle b_\nu(\Phi^0)\brangle=\Phi^0(\tau)\,.
$$
One finds
$$
\mu\partial_\tau{\bar I^0}(\tau,\mu)+\kappa({\bar I^0}(\tau,\mu)-\Phi^0)+\lambda({\bar I^0}(\tau,\mu)-\widetilde{\bar I^0}(\tau))=0\,,
$$
or, equivalently
$$
\begin{aligned}
{}&\mu\partial_\tau{\bar I^0}(\tau,\mu)+(\kappa+\lambda)({\bar I^0}(\tau,\mu)-\widetilde {\bar I^0}(\tau))=0\,,\qquad \tau>0\,,
\\
&{\bar I^0}(0,\mu)=\mu\blangle Q^0_\nu\brangle\,,\qquad\qquad\qquad\qquad\qquad\quad 0<\!\mu\!<1.
\end{aligned}
$$
Once $\widetilde {\bar I^0}$ is known, one recovers $I^0_\nu$ by using the semi-analytical formula \eqref{IntForm}, which in this case is
\begin{equation}\label{RT00}
\begin{aligned}
I^0_\nu(\tau ,\mu)=&\mathbf 1_{0<\mu<1}\mu e^{-\frac{(\kappa+\lambda)\tau }{\mu}}Q^0_\nu
+\mathbf 1_{0<\mu<1}\int_0^\tau e^{-\frac{(\kappa+\lambda)(\tau -y)}{\mu}}\frac{\kappa b_\nu(\widetilde {\bar I^0}(y))+\lambda\widetilde I^0_\nu(y)}{\mu}dy
\\
&+\mathbf 1_{-1<\mu<0}\int_\tau ^\infty e^{-\frac{(\kappa+\lambda)(y-\tau )}{|\mu|}}\frac{\kappa b_\nu(\widetilde {\bar I^0}(y))+\lambda\widetilde I^0_\nu(y)}{|\mu|}dy\,.
\end{aligned}
\end{equation}
This is precisely the solution of the Milne problem in the case of radiative transfer in a grey atmosphere.

In order to study the effect of the perturbation of the absorption defined above, we seek to compute
$$
I'_\nu(\tau ,\mu)=\frac{\partial I_\nu^\epsilon(\tau ,\mu)}{\partial\epsilon}\Big|_{\epsilon=0}\,,\quad\Phi'(\tau )=\frac{\partial\Phi^\epsilon(\tau )}{\partial\epsilon}\Big|_{\epsilon=0}\,,\quad u'(\tau ,\mu)=\frac{\partial u^\epsilon(\tau ,\mu)}{\partial\epsilon}\Big|_{\epsilon=0}\,.
$$
Let $\dot b_\nu(\Phi)=\partial_\Phi b_\nu(\Phi)$. Let us differentiate in $\epsilon$ both sides of the equations  \eqref{RTeps} at $\epsilon=0$: one easily finds that
$$
\begin{aligned}
\mu\partial_\tau I'_\nu+\kappa(I'_\nu-\dot b_\nu(\Phi^0)\Phi'_0)+\lambda(I'_\nu-\widetilde I'_\nu)+\delta\kappa_\nu(I^0_\nu-b_\nu(\Phi^0))=0\,,
\\
\kappa\blangle\widetilde I'_\nu\brangle-\kappa\underbrace{\blangle\dot b_\nu(\Phi^0)\brangle}_{=1}\Phi'+\blangle\delta\kappa_\nu(\widetilde I^0_\nu-b_\nu(\Phi^0))\brangle=0\,,
\end{aligned}
$$
and we can recast the second equality as
\begin{equation}\label{phiprime}
\Phi'(\tau )=\widetilde {\bar {I^0}}^\prime(\tau )+\tfrac1\kappa\blangle\delta\kappa_\nu(\widetilde I^0_\nu(\tau )-b_\nu(\Phi^0(\tau )))\brangle=\widetilde {\bar {I^0}}^\prime(\tau )+\tfrac1\kappa\blangle\delta\kappa_\nu(\widetilde I^0_\nu(\tau )-b_\nu(\widetilde {\bar I^0}(\tau )))\brangle\,.
\end{equation}

Averaging in frequency the RT equation, we find that
$$
\mu\partial_\tau {\bar {I^0}}^\prime(\tau ,\mu)+\kappa({\bar {I^0}}^\prime(\tau ,\mu)-\underbrace{\blangle\dot b_\nu(\Phi^0(\tau ))\brangle}_{=1}\Phi'(\tau ))+\blangle\delta\kappa_\nu(I^0_\nu(\tau ,\mu)-b_\nu(\Phi^0(\tau )))\brangle=0\,,
$$
and we can further eliminate $\Phi'(\tau )$ between the last two equations, to find,
\begin{equation}\label{uprime}
 \mu\partial_\tau {\bar {I^0}}^\prime(\tau ,\mu)+(\kappa+\lambda)({\bar {I^0}}^\prime(\tau ,\mu)-\widetilde {\bar {I^0}}^\prime(\tau ))=-\blangle\delta\kappa_\nu(I^0_\nu(\tau ,\mu)-\widetilde I^0_\nu(\tau ))\brangle\,.
\end{equation}
If $\delta\kappa$ is small everywhere but large near $\nu^*$ and if $\kappa+\lambda \ll 1$ so that the righthandside dominates the seccond term on the left, then \eqref{uprime} becomes
\[
 \mu\partial_\tau {\bar {I^0}}^\prime(\tau ,\mu)
 =-\delta\kappa(I^0_{\nu^*}(\tau ,\mu)-\widetilde I^0_{\nu^*}(\tau ))\,.
 \]  
Together with \eqref{phiprime} it gives the sign of $T'(z)$ by,
\begin{equation}\label{TempVar}
\Phi'(\tau )=
-\delta\kappa\Big|\int^\tau\widetilde{\frac{I^0_{\nu^*}(t ,\mu)}{\mu} }d t\Big|
+\frac{\delta\kappa}\kappa(\widetilde I^0_{\nu^*}(\tau )-b_{\nu^*}(\widetilde {\bar I^0}(\tau )))\,.
\end{equation}
when $\kappa\ll 1$ the second terms dominate the first one and it is possible to obtain a positive sign when $\widetilde I^0_{\nu^*}(\tau ) > b_{\nu^*}(\widetilde {\bar I^0}(\tau ))$.

\subsection{Numerical Simulations}
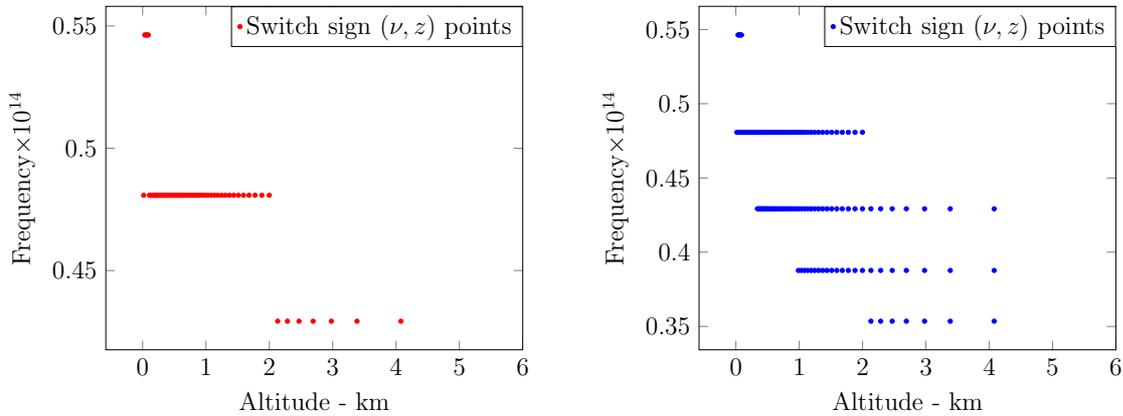
\begin{figure}[htbp]
\begin{minipage}[b]{0.45\textwidth}
\begin{center}
\begin{tikzpicture}[scale=0.8]
\begin{axis}[legend style={at={(1,1)},anchor=north east}, compat=1.3,
  xlabel= {Altitude - km},
  ylabel= {Frequency$\times 10^{14}$}
  , xmax=6]
%
\addplot[color=red,only marks, mark size=1pt] table [x index=0, y index=1]{fig/switch11.txt};
\addlegendentry{ Switch sign $(\nu,z)$ points}
%
\end{axis}
\end{tikzpicture}
\end{center}
\end{minipage}
\hskip0.5cm
\begin{minipage}[b]{0.45\textwidth}
\begin{center}
\begin{tikzpicture}[scale=0.8]
\begin{axis}[legend style={at={(1,1)},anchor=north east}, compat=1.3,
  xlabel= {Altitude - km},
  ylabel= {Frequency$\times 10^{14}$}
  , xmax=6]
%
\addplot[color=blue,only marks, mark size=1pt] table [x index=0, y index=1]{fig/switch12.txt};
\addlegendentry{ Switch sign $(\nu,z)$  points}
%
\end{axis}
\end{tikzpicture}
\end{center}
\end{minipage}
\caption{\label{fem5}Switch sign points: left when $\nu$ is increased at fixed $z$. Right : when altitude is increase at fixed $\nu$.  }
\end{figure}

We study the sign of \eqref{sign}.  As it is a complex surface of $\nu$ and $\tau$, we look at the frequencies and altitudes for which  \eqref{sign} changes sign when one of the two variable varies.
Figure \ref{fem5} shows the frequencies and altitudes where there is a change of sign, for the grey case $\kappa=0.5$.

Consider now $\kappa_\nu=0.5+\delta\kappa{\bf 1_{(\nu^1,\nu^2)}}$. We use an AD (automatic differentiation) extension of the computer program for \eqref{RT00} and display the sensitivity of $T$ with respect to $\delta\kappa$ for 3 ranges of support $(\nu^1,\nu^2)$. 

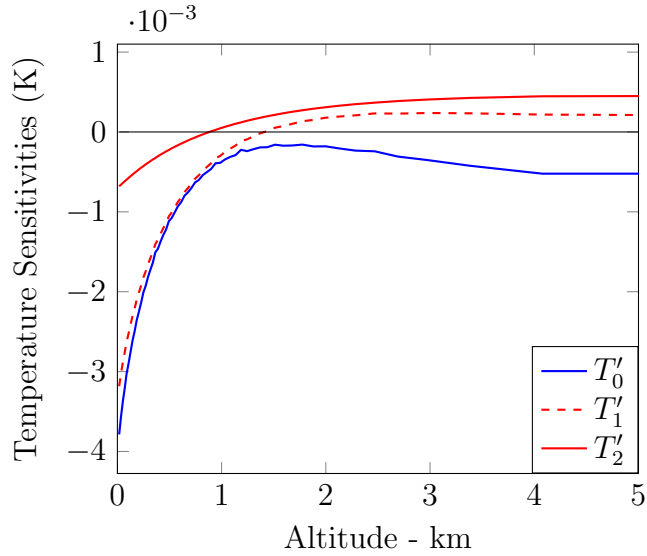
\begin{figure}[htbp]
\begin{center}
\begin{tikzpicture}[scale=1]
\begin{axis}[legend style={at={(1,0)},anchor=south east}, compat=1.3,
  xlabel= {Altitude - km},
  ylabel= {Temperature Sensitivities (K)}, xmin=0, xmax=5, ymax=0.0011]
\addplot[thick,solid,color=blue,mark=none,mark size=1pt] table [x index=0, y index=4]{fig/derivative.txt};
\addlegendentry{ $T'_0$}
\addplot[thick,dashed,color=red,mark=none,mark size=1pt] table [x index=0, y index=2]{fig/derivative.txt};
\addlegendentry{ $T'_1$}
\addplot[thick,solid,color=red,mark=none,mark size=1pt] table [x index=0, y index=3]{fig/derivative.txt};
\addlegendentry{ $T'_2$}
\addplot[thin,solid,color=black,mark=none,mark size=1pt] table [x index=0, y index=1]{fig/derivative.txt};

%
\end{axis}
\end{tikzpicture}
\end{center}
\caption{\label{fem6}  Scaled temperatures  sensitivities $z\to T'(z)$ versus altitude, computed with $\kappa_\nu=0.5+\delta\kappa{\bf 1_{(\nu^1,\nu^2)}}$ and $\delta\kappa$, at $\delta\kappa=0$ when case 0: $\nu^1=0.2$, $\nu^2=0.3$. Case 1: $\nu^1=0.3$, $\nu^2=0.4$. Case 2 is with $\nu^1=0.6$, $\nu^2=0.8$ . Notice the change of signs in the sensitivity $T'_2$.}
\end{figure}

On Figure \ref{fem6} it is seen that when $\nu^1=0.2$, $\nu^2=0.3$ the sensitivity of $T$ is negative, meaning that increasing $\kappa_\nu$ in the range $(\nu^1,\nu^2)$ triggers a decrease of temperature at all altitudes. 
Conversely with $\nu^1=0.6$, $\nu^2=0.8$ the sensitivity is positive above 800m meaning that an increase of $\kappa_\nu$ in that range triggers a temperature decrease below 800m and an increase above 800m.  The same is true for $T'_1$ but at 1500m instead of 800m.
%


\bibliographystyle{plain}
\bibliography{references}


\begin{lstlisting}[style=CStyle]

//  main.cpp : solves the radiative trnsfer equations
// for any absorption coefficient kappa
// does automatic differentiation if _NOAUTODIF is undefined
//  Created by Olivier Pironneau on 10/10/2021.
//

#include <iostream>
#include <fstream>
#include <cmath>
#include <string>
#include <time.h>
using namespace std;
#define sqr(x) (x*x)

#define _NOAUTODIF

#if defined(_NOAUTODIF)
    #define ddouble double
    #define ffabs(x) (fabs(x))
    const double dknu=0;
#else
    #include "ddouble.h"
    #define ffabs(x) (fabs(x.val[0]))
    const ddouble dknu(0.,1.); // means dknu=0.1 and AD computes d_T/d_dknu
#endif
const bool greycase=true;  // see beginning of main() for parameters. If false use below
const string kappa012="kappa0f";  // defines 3 runs with kappa vs nu in columns 2,3,4
const bool verbose=false;
const int Ntau=60; // nb points in tau
const int kmax=16;  // nb fixed point iterations
const double Z=1-exp(-12.); // max tau after change of var
const double SBsun =2.03e-5*(2*sqrt(2.))/0.7;  // scaled sunlight power
const double Tsun = 1.209;  // scaled sun temperature
const int jmaxmax=600; // max of max nb of points for integration in nu range
const int newton = 50; // to compute the temperature from int k*Botlzmann=int k*Imean
const double epsdycho=0.01, epsnewton=1.e-12;  // precision for dychotomy before Newton
const double tmin=1.e-10; // min t in ExpInt(t)
const double dtt = 0.005; //min integration step size in ExpInt integrals
const double kappamin=0.001;  // if kappa read is too small max it with kappamin
const int nt = 5; // min nb of integration step in anal formula
const double  ealb=0.3; //  and Earth albedo
const double alpha=1; // proportion of light affected at altitude 0 vs Z
const double ais=0., ars=0.; // max isotropic and Rayleigh scattering
const double tm1 = Z*0.6, tm2=Z*0.9;  // altitudes for scattering
double  nu[jmaxmax],// uneven discretization of [numin,numax]
        ai[jmaxmax], ar[jmaxmax], aux[jmaxmax];// isotropic and Rayleigh scattering
ddouble Inut[Ntau], // mu integral of I_nu
        Snut[Ntau],   // mu integral of mu^2 I_nu
        G[jmaxmax][Ntau], //  J(nu,tau)
        S[jmaxmax][Ntau], //  K(nu,tau)
        T[Ntau], kappanu[jmaxmax]; // T0 with kappanu0[]
string basedir("/Users/pironneau/Dropbox/afaire/Golse-Bardos/greenhouse4/");
string myresulttemperature("temperature");
string myresultmeanintensity("imean");
string mykappafile(basedir+kappa012+".txt"); // has nu, kappa0, kappa1, kappa2
int jmax;

struct ddouble2{ ddouble d1,d2;};

// don't use if kappa >18/Z
ddouble expint_E1(const ddouble t=1){
    // if your compiler has it or if you can link to gsl you may adapt this function
    // it computes E1(t)*B  when t<18
    double abst=ffabs(t);
    const int K = 9+(abst-1)*4;; // precision in the exponential integral function E1
    const double  gaNtaua =0.577215664901533; // special integration for log(t)
    if(abst<tmin) return 0.;// because integral_0^t(logx)dx ~0
    if(abst>18) {cout << "value of E_1 is incorrect with "<<t<<">18"<<endl; return 0;}
    ddouble ak=(t<0)?-t: t, soNtaue=-gaNtaua - log(abst)+ak;
    for(int k=2;k<K;k++){
        ak *= -abst*(k-1)/sqr(k);
        soNtaue += ak;
    }
    return soNtaue;
}
ddouble expint_E2(const ddouble t=1){
    ddouble aux = exp(-t) - t*expint_E1(t);
    return aux;
}
ddouble expint_E3(const ddouble t=1){
    return (exp(-t) - t*expint_E2(t))/2;
}
ddouble expint_E4(const ddouble t=1){
    return (exp(-t) - t*expint_E3(t))/3;
}
ddouble expint_E5(const ddouble t=1){
    return (exp(-t) - t*expint_E4(t))/4;
}

int readkappa(string mykappafile, int which=0){
    // on each line of kappafile: nu kappa0  & optional kappa1  kappa2
    ifstream kappafile(mykappafile);
    int j=-1;
    double kappaux=0, dummy, nuj=0.01;
    while((j++<=jmaxmax)&&(!kappafile.eof())){
        if(which==0) kappafile >> nuj >>kappaux >> dummy >> dummy >> dummy;
        else if(which==1) kappafile >> nuj >> dummy >>kappaux  >> dummy >> dummy;
        else if(which==2) kappafile >> nuj >> dummy >> dummy >>kappaux >> dummy;
        kappanu[j]  = fmax(kappaux,kappamin);
        nu[j]=nuj;
    }
    kappafile.close();
    jmax=j;
    return 0;
}

double Bsun(const double nu){ return SBsun*sqr(nu)*nu/(exp(nu/Tsun) -1);} // Boltzmann

ddouble BB(const double nu, const ddouble T){
    if(T<1.e-7) return 0;
    if(nu<1e-10) return T*sqr(nu);
    return sqr(nu)*nu/(exp(nu/T) -1);} // Boltzmann
ddouble dBB(const double nu, const ddouble T){
    if(T<1.e-7) return 0;
    if(nu<1.e-10) return sqr(nu);
    ddouble a = exp(nu/T); return a*sqr(nu*nu/(a -1)/T);
}

ddouble2 intBS(const int jnu,const double tau,const double tmin,const double tmax){
    // returns the convolution t-integral for the mean in mu of kappa*I_nu and mu^2 kappa*I_nu
    ddouble Imean=0, Imu2mean=0;
    ddouble kappa=kappanu[jnu];
    const double dt=fmin(dtt,nt/(tmax-tmin));
    for(double t=tmin;t<tmax;t+=dt){
        int it = int((Ntau-1)*t/Z); // parabolic length for frequency step
        double arnu4=ar[it]*sqr(nu[jnu]-0.8)*sqr(nu[jnu]-1.2)*(nu[jnu]>0.8)*(nu[jnu]<1.2)*40;
        ddouble  H0 = kappa*( BB(nu[jnu],T[it])*(1-arnu4) +
                        (ai[it]+1.125*arnu4)*G[jnu][it] -1.125*arnu4*S[jnu][it]) ;
        ddouble H2 = -0.375*arnu4*kappa*(G[jnu][it] -3*S[jnu][it]);
        if(kappa*(t-tau)!=0){
            Imean += dt*H0*(expint_E1(kappa*fabs(tau-t))+ealb*expint_E1(kappa*(tau+t)))/2;
            Imu2mean += dt*H0*(expint_E3(kappa*fabs(tau-t))+ealb*expint_E3(kappa*(tau+t)))/2;
            if(H2 !=0){
                Imean += dt*H2*(expint_E3(kappa*fabs(tau-t))+ealb*expint_E3(kappa*(tau+t)))/2;
                Imu2mean += dt*H2*(expint_E5(kappa*fabs(tau-t))+ealb*expint_E5(kappa*(tau+t)))/2;
            }
        }
    }
    ddouble2 mean; mean.d1=Imean; mean.d2 = Imu2mean;
    return mean;
}

ddouble thefunc(const ddouble rhs, const ddouble T0){
    ddouble myeq=-rhs;
    for(int j=1; j<jmax;j++){
         myeq += kappanu[j]*BB(nu[j],T0)*(nu[j]-nu[j-1]);
    }
    return myeq;
}
ddouble getTbydycho(const int i,const ddouble Tstart){
    ddouble T0=Tstart;
    if(Tstart<0.1) T0=0.1;
    ddouble Taux, myeq0=1,myeq1=-1, T1=T0, rhs=0;
    for(int j=1; j<jmax;j++){
         rhs+=kappanu[j]*G[j][i]*(nu[j]-nu[j-1]);
    }
    myeq0 =-rhs;myeq1=-rhs;
    while (myeq0>0){
        T0=T0/2;
        myeq0=thefunc(rhs,T0);
        if(verbose) cout<<T0<<" down "<<myeq0<<endl;
    }
    while (myeq1<0){
        T1=2*T1; myeq1=thefunc(rhs,T1);
        if(verbose) cout<<T1<<" up "<<myeq1<<endl;
    }
    while (T1-T0 > epsdycho){
        Taux=(T1+T0)/2;
        myeq0=thefunc(rhs,Taux);
        if(myeq0>0) T1=Taux; else T0=Taux;
        if(verbose) cout<<T0<<" middle "<<T1<<" "<<myeq0 << endl;
    }
    return (T1+T0)/2;
}

int genT(){
    for(int i=0;i<Ntau;i++){
        T[i] = getTbydycho(i,T[i]);
        ddouble presfunc = 1;
        int inewton =0;
        while(inewton++<newton && ffabs(presfunc) > epsnewton){
            ddouble T0 = T[i];
            ddouble left=0, rhs=0, deriv=0;
            double nu1=0;
            for(int j=1; j<jmax;j++){
                double dnu=nu[j]-nu[j-1]; // variable integral increment
                nu1=(nu[j]+nu[j-1])/2;
                rhs+=kappanu[j]*G[j][i]*dnu;
                left += kappanu[j]*BB(nu1,T0)*dnu;
                deriv += kappanu[j]*dBB(nu1,T0)*dnu;
            }
            presfunc = rhs-left;
            if(ffabs(deriv)>1e-10) T[i] = T0+presfunc/deriv;
            if(verbose) cout<<T[i]<<" newton "<<presfunc<<endl;
        }
        if(inewton>=newton)cout << "Newton precision doubtful "<<endl;
    }
    return 0;
}

int getT(){ // return temperature by Kirchhof's law when kappa is constant
    const double pi = 4*atan(1.);
    for(int i=0;i<jmax;i++){
        ddouble rhs=0;
        for(int j=1; j<jmax;j++)
            rhs+=G[j][i]*(nu[j]-nu[j-1]);
        T[i]=sqrt(sqrt(15*rhs))/pi;
    }
    return 0;
}

int getISnu(const int j){ // returns mean I_nu & mean mu^2 I_nu
    double nu1 = nu[j];
    ddouble kappa=kappanu[j];
    if(kappa<0.01) kappa=0.01;
    for(int i=0;i<Ntau;i++){
        double x=i*Z/(Ntau-1);
        ddouble2 aux = intBS(j,x,0,Z);
        Inut[i] = aux.d1 + Bsun(nu1)*(alpha*expint_E3(kappa*(x+Z)) + (1-alpha)*expint_E3(kappa*(Z-x))
                                      )/2;//              +ealb*alpha*expint_E3(kappa*(Z+x)))/2;
        Snut[i] = aux.d2 + Bsun(nu1)*(alpha*expint_E5(kappa*x) +(1-alpha)*expint_E5(kappa*(Z-x))
                                      )/2;//       +ealb*alpha*(expint_E5(kappa*(Z+x))))/2;
      }
    return 0;
}

int multiBlock(double initT)
{
    for(int i=0;i<Ntau;i++){ T[i]=initT;  // initialize
        for(int j=0; j<=jmax;j++){ G[j][i]=0; S[j][i]=0; }
        }
    for(int k=0;k<kmax; k++){  // fixed point loop: first update F
        for(int i=0;i<Ntau;i++){
              Inut[i]=0;  Snut[i]=0;
         }
        for(int j=0; j<jmax;j++){
            getISnu(j);
            for(int i=0;i<Ntau;i++){
                G[j][i]=Inut[i];
                S[j][i]=Snut[i];
            }
       }
  if( greycase )
#if defined(_NOAUTODIF)
    genT(); // replace by getT() to speed-up but know what you do!
#else
    genT();// Then update T.  AD need genT() and false with getT() because kappa varies
#endif
  else
    genT();
    ddouble normG=0, normS=0;
    for(int j=0; j<jmax;j++)
        for(int i=0;i<Ntau;i++){ normG+=ffabs(G[j][i]); normS+=ffabs(S[j][i]); }
    cout << "k= "<<k <<" "<<T[2]<<"  "<<normG<<"  "<<normS<<endl;
//      ofstream tempfile(basedir+"tempHistScatLow"+std::to_string(k)+".txt"); //study convergence
//        for(int i=1;i<Ntau;i++) tempfile<<-log(1-i*Z/(Ntau-1))<<"\t" <<T[i] <<endl;
    }
     return 0;
}

int main(int argc, const char * argv[]) {
    for(int it=0;it<Ntau;it++){   // defines scattering values vs tau
        double tm = it*Z/Ntau;
        ai[it] = ais*fmax(tm-tm1,0.)*fmax(tm2-tm,0.)*4/sqr(tm2-tm1);
        ar[it] = ars*fmax(tm-tm2,0.)/(Z-tm2);
    }
    for(int which=0;which<1;which++){ // change to which<1 for a single computation
    if(greycase){
        const double kappa0=0.5, numin = 0.05, numax = 15;
        jmax = 400;
        for(int j=0;j<jmax;j++){
            nu[j]= numax/(1+(jmax-j)*(numax-numin)/numin/jmax);  // uniform in wavelength
 //           kappanu[j]=kappa0;//+dknu*(nu[j]>0.2)*(nu[j]<0.3);
            kappanu[j]=kappa0*(nu[j]<6) + 0.1;
        }
    } else{
        readkappa(mykappafile,which);
//        for(int j=0;j<jmax;j++)  if((nu[j]<3./4)*(nu[j]>3./7)) kappanu[j]=fmax(0.5,1.5*kappanu[j]);
    }
    ofstream  resultfile =
        ofstream(basedir+myresulttemperature+to_string(which+10*(ais>0))+".txt");
    cout<<"\n iterations \t [T] near earth [T] far ||G|| and ||S||\n";
    double t0 = clock();
    multiBlock(0.);//,"tempmin");
    printf( " Time CPU = %10.6f\n", (clock() - t0)/CLOCKS_PER_SEC);
    cout<<"\n tau\t [T]:"<<endl;
    for(int i=1;i<Ntau;i++){
       cout << -log(1-i*Z/(Ntau-1))<<"\t"<<T[i] <<endl;
       resultfile<< -log(1-i*Z/(Ntau-1))<<"\t" // altitude
        <<T[i]<<endl;        // T(kappa) Milne by multigroup
    }
    ofstream  imeanz0 =
            ofstream(basedir+myresultmeanintensity+to_string(10*(ais>0)+which)+"0.txt");
        ofstream  imeanzZ =
                ofstream(basedir+myresultmeanintensity+to_string(10*(ais>0)+which)+"Z.txt");
        for(int j=10; j<jmax;j++){
            imeanz0<< 3/nu[j] <<" "<<1e5*G[j][1]<<endl; // mean intensity near the ground
            imeanzZ<< 3/nu[j] <<" "<<1e5*G[j][Ntau-1]<<endl;  // mean intensity at max altitude
        }
       resultfile.close();
        imeanz0.close();
    }
    return 0;
}

\end{lstlisting}
\end{document}